\documentclass[a4paper,11pt]{article}

\usepackage{graphicx}
\usepackage{latexsym}
\usepackage{amsmath}
\usepackage{amsthm}
\usepackage{amsfonts}
\usepackage{amssymb}
\usepackage{enumerate}

\usepackage[authoryear]{natbib}
  4  \oddsidemargin-20pt \evensidemargin-20pt
\allowdisplaybreaks[4]

\newcommand\N{\mathbb{N}}

\newcommand\R{\mathbb{R}}
\newcommand\Z{\mathbb{Z}}

\newcommand{\vect}[1]{\mathbf{#1}}  
\newcommand\Prob{\mathbb{P}}    
\newcommand\Exp{\mathbb{E}}     
\newcommand\ind{\mathbb{I}}     

\newcommand\Hc{\mathcal{H}}

\newcommand\Bb{\mathbb{B}}
\newcommand\Cb{\mathbb{C}}
\newcommand\Db{\mathbb{D}}
\newcommand\Eb{\mathbb{E}}
\newcommand\Gb{\mathbb{G}}

\DeclareMathOperator{\Cov}{Cov}

\DeclareMathOperator{\BL}{BL}

\newcommand{\linf}{\ell^\infty([0,1]^d)}

\newcommand\weak{\ \rightsquigarrow\ }

\newcommand\weakPi[1]{ \overset{\Prob}{\underset{ #1}{\,\rightsquigarrow}\,} }

\newcommand\Pconv{ \stackrel{ \Prob}{\rightarrow} }

\newtheoremstyle{normal}
{2ex}               
{3ex}               
{}                  
{}                  
{\bfseries} 
{}                  
{2pt}   
{\thmname{#1}\thmnumber{ #2.} \thmnote{(#3)}}

\newtheoremstyle{italic}
{2ex}
{3ex}
{\itshape}
{}
{\bfseries} 
{}
{2pt}
{\thmname{#1}\thmnumber{ #2.} \thmnote{(#3)}}

\theoremstyle{normal}
\newtheorem{definition}{Definition}[section]
\newtheorem{remark}[definition]{Remark}
\newtheorem{example}[definition]{Example}
\newtheorem{cond}[definition]{Condition}

\theoremstyle{italic}
\newtheorem{theorem}[definition]{Theorem}

\newtheorem{corollary}[definition]{Corollary}

\begin{document}

\title{Empirical and sequential empirical copula processes under serial dependence}

\author{Axel B\"ucher and  Stanislav Volgushev \\
Ruhr-Universit\"at Bochum \\
Fakult\"at f\"ur Mathematik \\
44780 Bochum, Germany \\
{\small e-mail: axel.buecher@ruhr-uni-bochum.de }\\
{\small e-mail: stanislav.volgushev@ruhr-uni-bochum.de}\\
}

\maketitle

\begin{abstract}
The empirical copula process plays a central role for statistical inference on copulas. Recently, \cite{segers2011} investigated the asymptotic behavior of this process under non-restrictive smoothness assumptions for the case of i.i.d. random variables. In the present paper we extend his main result to the case of serial dependent random variables by means of the powerful and elegant functional delta method.
Moreover, we utilize the functional delta method in order to obtain conditional consistency of certain bootstrap procedures.
 Finally, we extend the results to the more general sequential empirical copula process under serial dependence.
\end{abstract}


\section{Introduction}
\def\theequation{1.\arabic{equation}}
\setcounter{equation}{0}

Let $F$ be a $d$-variate continuous distribution function, $d\geq2$, with marginal distribution functions $F_p$ for $p\in\{1,\dots,d\}$. By the famous Theorem of Sklar [see \cite{sklar1959}] we can decompose $F$ as follows
\begin{align}	\label{sklar}
	F(\vect{x}) = C(F_1(x_1),\dots, F_d(x_d)), \quad \vect{x}=(x_1,\dots,x_d)\in\R^d,
\end{align}
where $C$ is the unique copula associated to $F$. By definition, $C$ is a $d$-variate distribution function on the unit cube $[0,1]^d$ whose univariate marginals are standard uniform distributions on the interval $[0,1]$. Equation \eqref{sklar} is usually interpreted in the way that the copula $C$ completely characterizes the information about the stochastic dependence contained in $F$. For an extensive exposition on the theory of copulas we refer the reader to the monograph \cite{nelsen2006}.

A fundamental tool for statistical inference on copulas is the so-called empirical copula. If $(\vect{X}_i)_{i\in\Z}$ forms a strictly stationary sequence of $d$-variate random vectors with copula $C$ and if $F_n$ denotes the empirical distribution function of $\vect{X}_1,\dots,\vect{X}_n$,  then the empirical copula is defined as
\begin{align*}
	C_n(\vect{u}) = F_n(F_{n1}^-(u_1), \dots, F_{nd}^-(u_d)), \quad \vect{u}\in[0,1]^d,
\end{align*}
where $F_{np}$ and $F_{np}^-$ denote the $p$-th marginal empirical cdf and its generalized inverse, respectively, for $p=1,\dots,d$. The asymptotic behavior of the corresponding empirical copula process
$$
	\Cb_n(\vect{u})=\sqrt{n}(C_n(\vect{u})-C(\vect{u}))
$$ 
was studied for $i.i.d.$ samples in \cite{gaenstut1987, ferradweg2004, tsukahara2005, vandwell2007, segers2011}, among others. Except for the results by the last named author, all previous derivations are based on the assumption that $C$ has continuous partial derivatives on the closed unit cube. \cite{segers2011} pointed out that there are hardly any popular copulas satisfying  this condition and showed that the results actually do hold under much less restrictive assumptions which are fulfilled for many copula models that are applied in practice. As we will demonstrate in Section \ref{sec:empcop}, the assumption made by \cite{segers2011} is in fact sufficient for weak convergence of the empirical copula process under various kinds of weak dependence.

In the case of weakly dependent strictly stationary samples $\vect{X}_1,\dots,\vect{X}_n$ only few results are available to date. To the best of our knowledge, \cite{douferlan2009} is the only reference where weak convergence of the empirical copula process is established for a sample of $\eta$-dependent random variables; again under the restrictive smoothness assumption of continuous partial derivatives on the whole unit cube discussed above. 
It is one of the main purposes of the present paper to extend the results of \cite{segers2011} and \cite{douferlan2009} to a weak convergence result for $\Cb_n$ under both a non-restrictive smoothness assumption and under various weak dependence concepts. At the same time, we provide simple conditions that allow for future generalizations of results on the empirical copula process. Our method of proof is based on the functional delta method [see also \cite{ferradweg2004, vandwell1996}], which implies that the weak convergence of the classical empirical process combined with a natural property of the limiting distribution and the mild smoothness assumption of \cite{segers2011} is already sufficient for the weak convergence of $\Cb_n$. As a by-product we rediscover the findings from \cite{segers2011} in a way that circumvents  sophisticated probabilistic arguments. Moreover, the functional delta method allows to prove consistency of bootstrap procedures for the empirical copula process that are constructed from consistent bootstrap approximations of the classical empirical process.

The second purpose of the present paper is an extension of the results of \cite{rueschendorf1976} for the sequential empirical process. 
More precisely, we consider weak convergence of the following two processes
\begin{align*}
	\Cb_n^\#(s,\vect{u}) &= \frac{1}{\sqrt{n}}\sum_{i=1}^{\lfloor sn \rfloor} \left( \ind\{ \vect{X}_i \leq F_n^-(\vect{u})  \} - C(\vect{u}) \right)   \\
		\Cb_n^+(s,\vect{u}) &=\frac{1}{\sqrt{n}}\sum_{i=1}^{\lfloor sn \rfloor} \left( \ind\{ \vect{X}_i \leq F_n^-(\vect{u})  \} - C_n(\vect{u}) \right),
\end{align*}
where $s\in[0,1]$ and $F_n^-(\vect{u})=(F_{n1}(u_1),\dots,F_{nd}^-(u_d))$. The first process is the one considered in \cite{rueschendorf1976} [except for a small error term of order $n^{-1}$ arising from the fact that $F_{np}^-(F_{np}(u_p))$ is not exactly equal to $u_p$] and it contains the usual empirical copula process as a special case, since $\Cb_n^\#(1,\vect{u})=\Cb_n(\vect{u})$. As in the case of the usual empirical copula process $\Cb_n$ described above, we both weaken the smoothness assumptions made in \cite{rueschendorf1976}  and we allow for various weak dependence concepts in order to obtain the weak convergence of $\Cb_n^\#$.

The second sequential empirical copula process $\Cb_n^+$ is for example applied in the context of testing for a constant dependence over time, see e.g. \cite{wiedehvanvog2011, remillard2010} or \cite{ruppert2011}. Except for the last named reference, the asymptotics of $\Cb_n^+$ are derived from the ones of $\Cb_n^\#$. However, the limiting process of $\Cb_n^+$ indicates that those assumptions should actually not be necessary. The latter author shows for the case of strictly stationary samples that $\Cb_n^+$ indeed converges weakly without \textit{any} smoothness assumption on $C$ at all.
In Section \ref{sec:seqempcop} we treat this problem by the functional delta method, which easily allows for more general weak dependence concepts. Note, that 
the methods  in \cite{wiedehvanvog2011} and \cite{remillard2010} are greatly improved in their applicability by our results.

The paper is organized as follows. In Section \ref{sec:empcop} we introduce the basic notation for the remaining part of the paper. Moreover, the usual empirical copula process, its weak convergence and a (blockwise) bootstrap approximation  under serial dependence is investigated. Some simple examples of applications of our results are given for illustration. In Section \ref{sec:seqempcop} we derive similar results for the sequential empirical copula processes.

\section{The empirical copula process}\label{sec:empcop}
\def\theequation{2.\arabic{equation}}
\setcounter{equation}{0}

\subsection{Weak convergence of the empirical copula process}

Let $(\vect{X}_i)_{i\in\Z}$ with $\vect{X_i}=(X_{i1},\dots, X_{id})$ be  a strictly stationary sequence of $d$-variate random vectors. Throughout the paper, the joint distribution function $F$ of $\vect{X_i}$ is assumed to have continuous marginal distributions $F_1, \dots, F_d$ and its copula is denoted by $C$. 
Further, let $U_{ip}=F_p(X_{ip})$ be the probability integral transform of the components of $\vect{X}_i$ and set $\vect{U_i}=(U_{i1},\dots, U_{id})$. Note that $\vect{U}_i$ is distributed according to $C$. Consider the empirical distribution functions
$$F_n(\vect{x})=n^{-1}\sum_{i=1}^n \ind\{ \vect{X}_i\leq \vect{x}\}, \quad G_n(\vect{u})=n^{-1} \sum_{i=1}^n \ind\{ \vect{U}_i \leq \vect{u} \}$$
for $\vect{x}\in\R^d,\vect{u}\in[0,1]^d$. For some distribution function $H$ on the real line let 
\begin{align*} H^{-}(p):=\begin{cases}  
				\inf\{x\in\R\,|\,H(x)\geq p \}, &\ 0<p\leq 1\\
				\sup\{x\in\R\,|\, H(x)=0\}, &\ p=0
				\end{cases} 
\end{align*}
denote the (left-continuous) generalized inverse function of $H$. The object of interest is the empirical copula defined as
$C_n(\vect{u}) = F_n(F_{n1}^-(u_1), \dots, F_{nd}^-(u_d))$, where  $F_{np}$ denotes the $p$-th marginal empirical cdf, $p=1,\dots,d$. 
For convenience, we will abbreviate the notation for the empirical copula by $C_n(\vect{u})= F_n(F_{n}^-(\vect{u}) )$, with $F_n^-(\vect{u})=(F_{n1}^-(u_1), \dots, F_{nd}^-(u_d))$.  Note that $G_n$ can not be computed from the data unless the marginal distributions $F_1,\dots,F_d$ are known. However, it can be shown that $C_n=G_n(G_{n1}^-,\dots,G_{nd}^-)$ [cf. \cite{segers2011}], which is the reason why we can prove the subsequent theoretical results by just considering the sequence $(\vect{U}_i)_{i\in\Z}$ and its empirical distribution function $G_n$. 

In order to derive the asymptotics of the corresponding empirical copula process
$\Cb_n=\sqrt{n}(C_n-C)$ we make the following assumption on the usual ($d$-variate) empirical process $\alpha_n=\sqrt{n}(G_n - C)$ based on $\vect{U}_1,\dots,\vect{U}_n$.

\begin{cond}\label{cond:weak}
	For the strictly stationary sequence $(\vect{X}_i)_{i\in\Z}$  it holds that 
	$\alpha_n=\sqrt{n}(G_n - C)$ weakly converges in $\linf$ to a tight, centered Gaussian field $\Bb_C$ concentrated on $\Db_0$, where
\begin{align*} 
	\Db_0=\left\{ \alpha\in C[0,1]^d\, | \, \alpha(1,\dots, 1)=0 \text{ and }  \alpha(\vect{x})=0 \, \text{ if some of the components of } \vect{x} \text{ are equal to } 0\right\}.
\end{align*}
\end{cond}

Here, and throughout the rest of this paper $\ell^\infty([0,1]^d)$ denotes the space of all uniformly bounded functions on the unit cube $[0,1]^d$, equipped with supremum norm $\|\cdot\|_\infty$.
Note that Condition \ref{cond:weak} is a standard result of empirical process theory for weakly dependent data, where the 
limiting process $\Bb_C$ is a tight, centered Gaussian process on $[0,1]^d$ with covariance 
\begin{align*}
	\Gamma(\vect{u},\vect{v})=\Cov(\Bb_C(\vect{u}),\Bb_C(\vect{v})) = \sum_{j\in\Z} \Cov(\ind\{\vect{U}_0\leq \vect{u}\}, \ind\{ \vect{U}_j\leq \vect{v}\} ).
\end{align*}
Sufficient conditions that entail Condition \ref{cond:weak}  have been proposed in terms of various weak dependence concepts:
\begin{itemize}
\item for $\alpha$-mixing (strongly mixing) conditions, see Theorem 2 in \cite{philpinz1980} or Theorem 7.3 in \cite{rio2000} for a more refined version. 
\item for $\beta$-mixing (absolutely regular) conditions, see e.g. Chapter 11.6 in \cite{kosorok2008} 
or \cite{doumasrio1995}.
\item for $\eta$-dependence, see Theorem 1 in \cite{douferlan2009}.
\item for multiple mixing properties see \cite{dehlduri2011} or \cite{duritusc2011}. 
\item for a similar result under long range dependence [at a different rate, though] see \cite{marinucci2005} and Remark \ref{rem:long} a) below.
\end{itemize}
It is the main purpose of this paper to prove weak convergence of $\Cb_n$ 
under the general Condition \ref{cond:weak} incorporating all the aforementioned citations. 

\begin{example}[ARMA$(p,q)$-processes] \label{ex:arma}
The $\R^d$-valued process $(\vect{X}_t)_{t\in\Z}$ is said to be an ARMA$(p,q)$-process if it is stationary and if it satisfies the equation
$$
	\sum_{i=0}^p B(i) \vect{X}_{t-i} = \sum_{j=0}^q A(j) \vect{Z}_{t-j}
$$
for every $t\in\Z$, where $B(i)\in \R^{d\times d}$ with $B(0)=Id$, $A(k)\in \R^{d\times r}$ $(r\geq 1)$ and where $\vect{Z}_t\in \R^r$  are i.i.d. centered random  vectors with covariance matrix $\Sigma$. For $z\in \Cb$ let $P(z)=\sum_{i=0}^p B(i)z^i$. It was shown in \cite{mokkadem1988} that if the absolute values of the zeros of the polyomial $\det P(z)$ are strictly greater than $1$ and if the probability law of $\vect{Z}_t$ is absolutely continuous with respect to the Lebesgue measure on $\R^r$, then $\vect{X}_t$ is $\beta$-mixing with $\beta(k)=O(e^{-\alpha k})$ for some $\alpha>0$. For instance by the results in \cite{kosorok2008}, this is sufficient for Condition \ref{cond:weak} to hold, with much to spare. 
\end{example}

For various further examples of stochastic processes  (e.g. GARCH processes, and various extensions thereof) satisfying any of the weak dependence concepts noted above we refer the reader to the corresponding citations.

In order to accomplish the step from the known weak convergence of $\alpha_n$ to the weak convergence of $\Cb_n$ some smoothness assumptions on $C$ are needed [even in the simple point-wise sense]. 
As pointed out by \cite{segers2011} the following smoothness condition on $C$ is nonrestrictive in the sense that it assures that the candidate limiting process of $\Cb_n$ has continuous trajectories. 

\begin{cond} \label{cond:pd}
For $p=1\dots,d$ the partial derivatives $\partial_p C(\vect{u})$  exist and are continuous on  $U_p=\{\vect{u} \in [0,1]^d\,|\,u_p \in (0,1)\}.$
\end{cond}

Now consider the mapping 
\begin{align*} 
\Phi: \begin{cases}
	\Db_\Phi &\rightarrow \ \ \ell^\infty[0,1]^d\\
	H &\mapsto \ \ H(H_1^-,\dots,H_d^-)=H(H^-),
\end{cases}
\end{align*}
where $\Db_\Phi$ denotes the set of all distribution functions $H$ on $[0,1]^d$ whose 
marginal cdfs $H_p$ satisfy $H_p(0)=0$. Observe that 
\[
	\Cb_n=\sqrt{n}(C_n-C)=\sqrt{n}(G_n(G_{n1}^-,\dots,G_{nd}^-)-C)=\sqrt{n}(\Phi(G_n)- \Phi(C)).
\]
Therefore, in order to derive the weak convergence of $\Cb_n$ it suffices to establish Hadamard-differentiability of $\Phi$ tangentially to suitable subspaces and to invoke the functional delta method. This Hadamard-differentiability  has been investigated in \cite{ferradweg2004} under the restrictive assumption that $C$ has continuous partial derivatives on the whole unit cube $[0,1]^d$ [cf. \cite{vandwell1996}, page 389]. The following Theorem relaxes this assumption to Condition \ref{cond:pd}, its proof is given at the end of this section.

\begin{theorem}\label{theo:C_hadamard}
Suppose Condition \ref{cond:pd} holds.
Then $\Phi$ is Hadamard-differentiable at $C$ tangentially to $\Db_0$.
Its derivative at $C$ in $\alpha\in \Db_0$ is given by
\begin{align*}
\big(\Phi'_C(\alpha)\big)(\vect{u}) = \alpha(\vect{u}) - \sum_{p=1}^d \partial_pC(\vect{u})\ \alpha(\vect{u}^{(p)}),
\end{align*}
where $\partial_p\,C$ is defined as 0 on the set $\{\vect{u}\in[0,1]^d\, |\, u_p\in\{0,1\} \}$ and where $\vect{u}^{(p)}=(1,\dots,1,u_p,1,\dots,1)$ denotes the vector, where all entries of $\vect{u}$ except the $p$-th are replaced by $1$.
\end{theorem}

The functional delta method, see e.g. Theorem 3.9.4 in \cite{vandwell1996}, immediately yields the following result. It has been stated in \cite{segers2011} for the i.i.d. case and in \cite{douferlan2009} for the serial dependent case, where the latter authors impose restrictive smoothness conditions on $C$ [and investigate the Skorohod-space $D([0,1]^d)$, which implies weak convergence with respect to the $\sup$-norm if the limiting process has a version that has continuous trajectories almost surely]. 

\begin{corollary} \label{corr:empcop}
Suppose Conditions \ref{cond:weak} and \ref{cond:pd} hold. Then the empirical copula process $\Cb_n=\sqrt{n}(C_n-C)$  weakly converges towards a Gaussian field $\Gb_C$,
\begin{align*} \Cb_n=\sqrt{n}({C}_n-C) ~\weak~ \mathbb{G}_C ~~~\text{ in }~~~ l^{\infty}([0,1]^d), \end{align*}
which may be expressed as
\begin{align} \label{G_C}
\mathbb{G}_C(\vect{u})=\mathbb{B}_C(\vect{u}) - \sum_{p=1}^d \partial_pC(\vect{u})\ \mathbb{B}_C(\vect{u}^{(p)}).
\end{align} 
\end{corollary}

\begin{remark}\label{rem:long}~\\
a) 
More generally, the functional delta method allows the derivation of weak convergence results under different rates of convergence of the empirical process. For instance, in the case of long range dependent stationary sequences, it was shown by \cite{marinucci2005} in the bivariate case that $r_n(G_n-C)$ weakly converges in $D([0,1]^2)$ to some limit process $\tilde \Bb_C\in \Db_0$, where the rate of convergence $1/r_n\rightarrow0$ satisfies $r_n=o(n^{1/2})$. If future research reveals more general results on the usual empirical process than the one by \cite{marinucci2005} our method of proof easily extends these results to the empirical copula process.

\medskip 

b) Note that a direct extension of the method of proof employed by \cite{segers2011} to the case of weakly dependent stationary sequences would need additional probabilistic arguments that differ depending on the underlying data structure. On the other hand, the functional delta method allows to separate the probabilistic part (weak convergence of the empirical distribution function and properties of the limiting process) and the analytical part (compact differentiability) of the problem hereby considerably simplifying the proof of the weak convergence of $\Cb_n$.
\end{remark}

The following example is a simple, but typical application of Corollary \ref{corr:empcop}.

\begin{example}\label{ex:spearman}\textbf{(Asymptotics of Spearman's rho)}\\
The multivariate population version of Spearman's rho is defined as
$$\rho=\frac{d+1}{2^d-(d+1)}\left( 2^d \int_{[0,1]^d} C(\vect{u})\, d\vect{u} - 1 \right),$$
see \cite{schmschm2007}, and may be estimated by the simple plug-in nonparametric estimator
$$\rho_n=\frac{d+1}{2^d-(d+1)}\left( 2^d \int_{[0,1]^d} C_n(\vect{u})\, d\vect{u} - 1 \right).$$
Observing continuity of the mapping $C\mapsto \rho(C)$ we immediately obtain the result that
$$\sqrt{n}(\rho_n-\rho) \weak N(0,\sigma^2), \qquad \sigma^2=\left(\frac{d+1}{2^d-(d+1)}\right)^{2}2^{2d}\int_{[0,1]^d}\int_{[0,1]^d} \Eb\left[\Gb_C(\vect{u})\Gb_C(\vect{v}) \right]\,d\vect{u}\,d\vect{v}.$$
\end{example}

\subsection{Bootstrapping the empirical copula process}
In the present section we apply the functional delta method in order to obtain consistent bootstrap procedures for the empirical copula process under a general assumption on a given bootstrap procedure for the usual $d$-variate empirical process. As an illustration of our results, we provide consistency of a simple (blockwise) bootstrap approximation of the limiting field $\Gb_C$ in the setting of $\alpha$-mixing data. More precisely, let $M_1,\dots, M_{n}$ be non-negative, real-valued random variables with $\sum_{i=1}^n M_i =n$ and define
\[
	F_{n,b}(\vect{x})=n^{-1}\sum_{i=1}^n M_i\ \ind\{ \vect{X}_j \leq \vect{x}\}, \qquad G_{n,b}(\vect{u})=n^{-1} \sum_{i=1}^n M_i\ \ind\{ \vect{U}_j \leq \vect{u}\}.
\]
For instance, if $(M_1,\dots,M_n)$ follows a multinomial distribution with parameters $n,n^{-1},\dots,n^{-1}$, we obtain the usual resampling bootstrap of an i.i.d. sample, i.e. the bootstrap based on sampling from the data with replacement. On the other hand, the construction above also allows to handle bootstrap procedures such as the blockwise bootstrap considered in Example \ref{ex:boot2} below.

In order to construct a consistent bootstrap procedure for the empirical copula process, we need the following condition.  

\begin{cond}\label{cond:bootweak}
Additionally to the assumptions made in Condition \ref{cond:weak}, suppose that $\alpha_{n,b}=\sqrt{n}(G_{n,b}-G_n)$ weakly converges in $\linf$ to $\Bb_C$ conditional on the data in probability, notationally $\alpha_{n,b}\weakPi{M} \Bb_C$.
\end{cond}

Weak convergence conditional on the data in probability ($\weakPi{M}$-convergence) is understood in the Hoffmann-J\o rgensen sense as defined in \cite{kosorok2008}, that is $\alpha_{n,b}\weakPi{M} \alpha$ if and only if 
\begin{enumerate}[(i)]
\item  $\sup_{f\in \BL_1} \left| \Exp_M f(\alpha_{n,b}) - \Exp f(\alpha) \right| \rightarrow 0$ in outer probability, 
\item $\Exp_M f(\alpha_{n,b})^* - \Exp_M f(\alpha_{n,b})_* \Pconv 0$ for all $f\in \BL_1$,
\end{enumerate}
where $\BL_1$ denotes the set of all Lipschitz-continuous functions $f:\linf\rightarrow\R$ that are uniformly bounded by $1$ and have Lipschitz constants bounded by $1$, and where the asterisks in (ii) denote measurable majorants (and minorants, respectively) with respect to the joint data $(\vect{X}_1,\dots,\vect{X}_n,M_1,\dots,M_n)$. We give three examples under which Condition \ref{cond:bootweak} is met.

\begin{example}\label{ex:boot1} \textbf{(Bootstrap in the i.i.d. case)} 
\begin{enumerate}[a)]
	\item Suppose $\vect{X}_1,\dots,\vect{X}_n$ are i.i.d. Moreover, assume that $\xi_1,\dots,\xi_n$ are i.i.d., non-negative, real-valued random variables that are independent of the $\vect{X}_i$ and satisfy the conditions $\Exp \xi_1=1, \Exp \xi_1^2 =1$ and $\|\xi_1\|_{2,1}=\int_0^\infty \sqrt{\Prob(|\xi_1|>x)}\,dx<\infty$. Define the random variables $M_1,\dots,M_n$ as $M_i := n\xi_i/(\sum_{i=1}^n \xi_i), i=1,\dots,n$. Then Theorem 2.6 in \cite{kosorok2008} implies that Condition \ref{cond:bootweak} holds.
	\item Assume that $\vect{X}_1,\dots,\vect{X}_n$ are i.i.d. and that the vector $(M_1,\dots,M_n)$ follows a multinomial distribution with parameters $n,n^{-1},\dots,n^{-1}$. Again, by Theorem 2.6 in \cite{kosorok2008} Condition \ref{cond:bootweak} holds.
\end{enumerate}	
\end{example}

\begin{example}\label{ex:boot2}	\textbf{(Bootstrap in the strongly mixing case)}\\
Suppose $(\vect{X}_i)_{i\in\Z}$ is a strictly stationary, strongly mixing sequence whose $\alpha$-mixing coefficients satisfy the summability condition $\sum_{i=0}^\infty(i+1)^{16(d+1)}\alpha^{1/2}(i)<\infty$ as considered in \cite{buehlmann1993}.
	In this setting, the usual bootstrap based on sampling with replacement from the data (see Example \ref{ex:boot1} b)) is not able to capture the serial dependence of the sequence. The blockwise bootstrap  considered in \cite{buehlmann1993} provides one possible solution. It is based on the following idea. Instead of sampling (with replacement) single data points, we consider $n-l+1$ blocks of length $l$ which have the form $(\vect{X}_1,\dots,\vect{X}_l), (\vect{X}_2,\dots,\vect{X}_{l+1}),\dots (\vect{X}_{n-l+1},\dots,\vect{X}_n)$. The bootstrap data are formed by sampling with replacement $k=n/l$ of those blocks. This works for $k$ being an integer, if that is not the case we simply make the last selected block smaller. If $l=l_n$ is chosen appropriately, this procedure is able to capture the serial dependence structure in the data. More formally, we define the following bootstrapped empirical distribution functions [for the sake of a clear exposition we only consider the case $k=n/l\in\N$]
\[
	F_{n,b}(\vect{x})=k^{-1}\sum_{i=1}^{k} l^{-1} \sum_{j=S_i+1}^{S_i+l} \ind\{ \vect{X}_j \leq \vect{x}\}, \qquad G_{n,b}(\vect{u})=k^{-1}\sum_{i=1}^{k} l^{-1} \sum_{j=S_i+1}^{S_i+l} \ind\{ \vect{U}_j \leq \vect{u}\},
\]
where $n=kl, l=l_n=o(n), l_n\rightarrow\infty$ and $(S_i)_{i\in\N}$ are i.i.d. uniformly distributed on the set $\{0,\dots,n-l\}$. Note that in the case $n = kl$ the quantity $F_{n,b}(\vect{x})$ can also be represented as $F_{n,b}(\vect{x}) = n^{-1}\sum_{i=1}^n M_i\ \ind\{ \vect{X}_j \leq \vect{x}\}$ where $M_j := \# \{i = 1,\dots,k\, |\, j-l \leq S_i \leq j-1 \}.$ Thus we are in the setting described at the beginning of this section.
Given that $l_n \rightarrow \infty,l_n=O(n^{1/2-\varepsilon})$ for some $\varepsilon>0$, Corollary 3.1 in \cite{buehlmann1993} states that $\alpha_{n,b}\weak \Bb_C$ almost surely in $(D([0,1]^d),\|\cdot\|_\infty)$, where the latter space is equipped with the $\sigma$-field induced by the open balls. Note that $\alpha_{n,b}$ is measurable with respect to the ball-$\sigma$-field and the theory for weak convergence on metric spaces spaces equipped with non-Borel sigma fields [see e.g. \cite{pollard1984}, Chapter IV] applies. Thus the delta method for weak convergence in the Hoffmann-J\o rgensen sense is not directly applicable. However, a closer look at Theorem 1.7.2 and the proof of Theorem 2.9.6 in \cite{vandwell1996} reveals that we can enlarge both the space $D([0,1]^d)$ to $\linf$ and the $\sigma$-field to the Borel $\sigma$-field and obtain $\alpha_{n,b}=\sqrt{n}(G_{n,b} - G_n) \weakPi{S} \Bb_C$ in $\linf$ in the Hoffmann-J\o rgensen sense.
\end{example}

Now, defining 
\[ 
	C_{n,b}(\vect{u})=F_{n,b}(F_{n,b}^{-}(\vect{u}))=G_{n,b}(G_{n,b}^{-}(\vect{u}))
\]
as the bootstrap empirical copula, the functional delta method for the bootstrap, see Theorem 12.1 in \cite{kosorok2008}, yields the following result.

\begin{corollary}\label{corr:bootweak}
	Under condition  \ref{cond:pd} and \ref{cond:bootweak} it holds 
	\[ 
		\Cb_{n,b}=\sqrt{n}(C_{n,b} - C_n) \weakPi{M} \Gb_C~~~\text{ in }~~~ \linf.
	\]
\end{corollary}

This result has also been stated in \cite{buecdett2010} under the special setting of Example \ref{ex:boot1} a) and in \cite{ferradweg2004} under the setting of Example \ref{ex:boot1} b), where both citations use the restrictive smoothness assumption that the copula has continuous partial derivatives on the whole unit cube. 

\begin{remark}
In the i.i.d. case the literature provides an alternative multiplier bootstrap procedure based on the estimation of the partial derivatives of the copula, see e.g. \cite{remiscai2009, segers2011}. Consistency is usually stated in a weaker unconditional sense [see \cite{buecher2011} for an extension to the conditional Hoffmann-J\o rgensen sense]. The investigation of this bootstrap procedure in the serial dependent case is beyond the scope of the present paper since it cannot be handled by the functional delta method. We refer the reader to \cite{ruppert2011} for the consideration of the strongly mixing case.
\end{remark}

Finally, we consider one of many possible applications of Corollary \ref{corr:bootweak}.

\begin{example}\label{ex:appl} \textbf{(Testing for symmetry of a copula)} \\ 
Suppose one is interested in a test for the (bivariate) hypothesis
$\Hc_0: C(u,v)=C(v,u)$ for all $u,v\in[0,1]$. Recently, \cite{gennesque2011} proposed to base a test for $\Hc_0$ by considering functionals of the following process
$$\Db_n(u,v)=\sqrt{n} (C_n(u,v) - C_n(v,u))\stackrel{\Hc_0}{=}\Cb_n(u,v)-\Cb_n(v,u),$$
which weakly converges to  $\Db_C(u,v)= \Gb_C(u,v)-\Gb_C(v,u)$
in $\ell^\infty([0,1]^2)$ by the continuous mapping theorem. By the bootstrap continuous mapping theorem, see \cite{kosorok2008}, we obtain
$$\Db_{n,b}(u,v)=\sqrt{n} (C_{n,b}(u,v) - C_{n,b}(v,u))\weakPi{S} \Db_C(u,v)$$
independently of the hypothesis. One test statistic investigated by \cite{gennesque2011} is defined as $T_n=\|\Db_{n}\|_\infty$, and we reject $\Hc_0$ if $T_n$ is larger than the $1-\alpha$ quantile of the law of $T_C=\|\Db_C\|_\infty$, which depends on the unknown copula in a complicated manner. Approximate quantiles can be obtained by a bootstrap sample of $T_{n,b}=\|\Db_{n,b}\|_\infty$.
\end{example}

\subsection{Proof of Theorem \ref{theo:C_hadamard}. }
Let $\Eb$ denote the set of distribution functions $F$ on $[0,1]$ with $F(0)=0$ and define $\Eb^-$ as the set of all generalized inverse functions $F^-$ with $F\in\Eb$. Now decompose 
 $\Phi=\Phi_3\circ\Phi_2\circ\Phi_1$, where
\begin{align} \label{Phi123}
\Phi_1&: \begin{cases}
	\Db_\Phi &\rightarrow \ \ \Db_{\Phi_2}=\Db_\Phi\times\Eb^d \\
	H &\mapsto \ \ (H, H_1, \dots, H_d)
\end{cases}\nonumber \\
\Phi_2&: \begin{cases}
	\Db_{\Phi_2} &\rightarrow \ \ \Db_{\Phi_3} = \Db_\Phi\times(\Eb^-)^d   \\
	(H,F_1,\dots,F_d) &\mapsto \ \ (H, F_1^-,\dots F_d^-) 
\end{cases} \\
\Phi_3&: \begin{cases}
	\Db_{\Phi_3} &\rightarrow \ \ \linf   \\
	(H,G_1,\dots,G_d) &\mapsto \ \ H\circ(G_1,\dots, G_d).
\end{cases}\nonumber
\end{align}
The first mapping $\Phi_1$ is Hadamard-differentiable at $C$ with derivative $\Phi_{1,C}'=\Phi_1$ since it is linear and continuous. 

Considering the second mapping let $U\in\Eb$  be the identity, then $\Lambda: F\mapsto F^-$ with $F\in\Eb$  is Hadamard-differentiable at $U$ tangentially to the set $\Eb_{0}=\{ F \in C[0,1]\ : \ F(0)=F(1)=0\}$  with derivative $\Lambda'_U(h)=-h$. To see this let $t_n\rightarrow0$, $h_n\in \ell^\infty([0,1])$ with $h_n\rightarrow h\in \Eb_{0}$ and $U+t_nh_n\in \Eb$. We have to show that 
\begin{align}\label{PsiU}
\sup_{x\in[0,1]} \left| \frac{\Lambda(U+t_nh_n)-\Lambda(U) }{t_n} - \Lambda'_U(h) \right|(x) = \sup_{x\in[0,1]} \left| \frac{(U+t_nh_n)^-(x) - x }{t_n} + h(x)\right| \rightarrow 0.
\end{align}

We deal with the case $x=0$ first and show that
\begin{align*}
A_n=(U+t_nh_n)^-(0)=\sup\{x\in[0,1] \ : \ x+t_nh_n(x)=0\} =o(t_n).
\end{align*} 
Suppose $A_n$ is non-zero, eventually. Choose a sequence $x_n\in[A_n/2,A_n]$ satisfying $x_n+t_nh_n(x_n)=0$, i.e. $h_n(x_n)=-x_n/t_n$. Since $h_n$ is bounded we obtain $x_n=O(t_n)=o(1)$. By uniform convergence of $h_n$ and continuity of $h$ this implies $-x_n/t_n=h_n(x_n)\rightarrow h(0)=0$, that is $x_n=o(t_n)$ and therefore $A_n=o(t_n)$.

Now consider the case $x>0$. Since $U+t_nh_n$ is a distribution function with $(U+t_nh_n)(0)=0$ it follows that $\xi_n(x)=(U+t_nh_n)^{-}(x)\in(0,1]$ for all $x\in(0,1]$. Set $\varepsilon(x)=t_n^2\wedge \xi_n(x)>0$, then
\begin{align*}
(U+t_nh_n)(\xi_n(x)-\varepsilon(x)) \leq x \leq (U+t_nh_n)(\xi_n(x))
\end{align*}  
for all $x\in(0,1]$ by the definition of the generalized inverse function. This implies
\begin{align*}
-t_nh_n(\xi_n(x))\leq \xi_n(x)-x \leq -t_nh_n(\xi_n(x)-\varepsilon(x)) + t_n^2.
\end{align*}
Since $h_n$ is uniformly bounded we obtain $\xi_n(x)\rightarrow x$ uniformly in $x\in(0,1]$. Divide the last equation by $t_n$ and use uniform convergence of $h_n$  and continuity of $h$ to conclude that $(\xi_n(x)-x)/t_n\rightarrow -h(x)$ uniformly in $x\in(0,1]$. Together with the case $x=0$ this yields \eqref{PsiU}.

It now easily follows that $\Phi_2$ is Hadamard-differentiable at $(C,U,\dots,U)$ tangentially to $\Db_0\times \Eb_{0}^d$ with derivative 
\begin{align*}
\Phi'_{2,(C,U,\dots,U)}(\alpha,f_1,\dots,f_d)=(\alpha,-f_1,\dots,-f_d).
\end{align*}

Last but not least consider $\Phi_3$. We assert that $\Phi_3$ is Hadamard-differentiable at $(C,U,\dots,U)$ tangentially to $\Db_0\times \Eb_0^d$
with derivative 
\begin{align*}
	\Phi'_{3,(C,U,\dots,U)}(\alpha,g_1,\dots,g_d)(\vect{x}) = \alpha(\vect{x})+\sum_{p=1}^d \partial_p C(\vect{x})g_p(x_p).
\end{align*}
To see this let $t_n\rightarrow0$ and $(\alpha_n, g_{n1}, \dots, g_{nd} )\in \linf \times (\ell^\infty([0,1]))^d$ with $(\alpha_n, g_{n1},\dots, g_{nd}) \rightarrow (\alpha,g_1,\dots,g_d)\in \Db_0\times \Eb_0^d$ such that $(C+t_n\alpha_n, U+t_ng_{n1},\dots, U+t_ng_{nd}) \in \Db_\Phi\times (\Eb^-)^d$. Now decompose
\begin{align*}
t_n^{-1} \{ \Phi_3(C+t_n\alpha_n, U+t_ng_{n1},\dots, U+t_ng_{nd}) - \Phi_3(C,U,\dots,U) \} = L_{n1}+L_{n2}
\end{align*}
where 
\begin{align*}
L_{n1} &= t_n^{-1} \{ C \circ (U+t_ng_{n1},\dots, U+t_ng_{nd}) - C \} \\
L_{n2} &= \alpha_n\circ (U+t_ng_{n1},\dots, U+t_ng_{nd})
\end{align*}
and consider both terms separately. Exploiting the facts that $\alpha_n$ and $(g_{n1},\dots,g_{nd})$ converge uniformly and that $\alpha$ is uniformly continuous one can conclude that $||L_{n2}-\alpha||_\infty = o(1)$. Concerning $L_{n1}$ we have to deal with different cases regarding the number of coordinates on the boundary of $[0,1]^d$. If $\vect{x}\in(0,1)^d$ a Taylor expansion of $C$ at $\vect{x}$ yields 
\begin{align*}
L_{n1}(\vect{x})=\sum_{p=1}^d \partial_pC(\vect{x}) g_{np}(x_p)  + r_n(\vect{x}), 
\end{align*}
where the error term $r_n$ can be written as
\begin{align*}
	r_n(\vect{x})=\sum_{p=1}^d \big(\partial_pC(\vect{y}) - \partial_pC(\vect{x})\big) g_{np}(x_p) 
\end{align*}
with some intermediate point $\vect{y}=\vect{y}(n)$ between $\vect{x}$ and $(U+t_ng_{n1},\dots,U+t_ng_{nd})(\vect{x})$. The main term uniformly converges to $\sum_{p=1}^d \partial_pC(\vect{x}) g_{p}(x_p)$ as asserted, hence it remains to show that the error term converges to $0$ uniformly in $\vect{x}$. To see this, let $\varepsilon>0$. Using uniform convergence of $g_{np}$, uniform continuity of $g_p$ and the fact that $g_p(0)=g_p(1)=0$ one can conclude that there exists $\delta>0$ such that $|g_{np}(x_p)| \leq \varepsilon/2$ for all $x_p<\delta$ and $x_p>1-\delta$ and all $p=1,\dots,d$. Since partial derivatives of copulas are bounded by $1$ the $p$-th summand of the error term is uniformly bounded by $\varepsilon$ for $x_p<\delta$ and $x_p>1-\delta$.  On the set $[0,1]^d\setminus\{x_p<\delta\text{ or } x_p>1-\delta\}$ the partial derivative $\partial_pC$ is uniformly continuous which yields the desired convergence under consideration of  $\vect{y}(n)\rightarrow\vect{x}$ and boundedness of $g_p$. Since $\varepsilon>0$ was arbitrary the case $\vect{x}\in(0,1)^d$ is finished.

We introduce the notation $\vect{x}_{(p)}=(x_1,\dots,x_{p-1},1,x_{p+1},\dots,x_d)$ and now consider the case $\vect{x}=\vect{x}_{(p)}$ with $x_q\in(0,1)$ for $q\ne p$. Decompose
$L_{n1}(\vect{x}) = L_{n1}^{(1)}(\vect{x})+L_{n1}^{(2)}(\vect{x})$
where
\begin{align*}
L_{n1}^{(1)}(\vect{x})&= t_n^{-1} \{ C( (x_1+t_ng_{n1}(x_1),\dots,x_d+t_ng_{nd}(x_d))_{(p)} ) - C(\vect{x}_{(p)})\} \\
L_{n1}^{(2)}(\vect{x})&= t_n^{-1} \{ C( x_1+t_ng_{n1}(x_1),\dots,x_d+t_ng_{nd}(x_d) ) - C( (x_1+t_ng_{n1}(x_1),\dots,x_d+t_ng_{nd}(x_d))_{(p)} ) \}.
\end{align*}
By the same arguments as in the previous case a Taylor expansion in $\vect{x}_{(p)}$
and a discussion of the remainder yields
\begin{align*}
L_{n1}^{(1)}(\vect{x})= \sum_{q\ne p} \partial_qC(\vect{x})g_q(x_q) + o(1) = \sum_{q=1}^d \partial_qC(\vect{x}) g_q(x_q)  + o(1)
\end{align*} 
uniformly in $\vect{x}$ (note that  $g_p(1)=0$). Lipschitz-continuity of $C$ implies $|L_{n1}^{(2)}(\vect{x}) | \leq |g_{np}(1)| \rightarrow |g_p(1)| = 0$ and the second case is finished. If more than one coordinate is $1$ while the others stay in $(0,1)$ the argumentation is analogous.

If one (or more) of the coordinates of $\vect{x}$ equals $0$, say $x_p=0$ use $C(\vect{x})\equiv0$  and Lipschitz-continuity of $C$ to estimate 
\begin{align*}
|L_{n1}(\vect{x}) |  \leq |g_{np}(0)| \rightarrow 0.
\end{align*}
If $x_p=0$ we have $\partial_qC(\vect{x})=0$ for $q\ne p$ and $g_p(x_p)=g_p(0)=0$, which entails
$0=\sum_{p=1}^d \partial_pC(\vect{x}) g_{p}(x_p)$ and thus the assertion.

The remaining case $\vect{x}=(1,\dots,1)$ again follows by Lipschitz-continuity, the details are omitted.
To conclude, $\Phi_3$ is Hadamard-differentiable as asserted.

Now apply the chain rule, Lemma 3.9.3 in \cite{vandwell1996}, to $\Phi=\Phi_3\circ\Phi_2\circ\Phi_1$ to conclude the assertion of the Lemma. \qed

\section{The sequential empirical copula process}\label{sec:seqempcop}
\def\theequation{3.\arabic{equation}}
\setcounter{equation}{0}

Under the setting of Section \ref{sec:empcop} consider the following partial sum  
$$G_n^{\#}(s,\vect{u})=n^{-1} \sum_{i=1}^{\lfloor sn\rfloor} \ind\{ \vect{U}_i \leq \vect{u} \}$$
for $s\in[0,1]$ and  $\vect{u}\in[0,1]^d$. The process 
$$\alpha_n^{\#}(s,\vect{u})=n^{-1/2}\sum_{i=1}^{\lfloor sn \rfloor} \left( \ind\{ \vect{U}_i \leq \vect{u} \} - C(\vect{u}) \right) = \sqrt{n}\left( G_n^{\#}(s,\vect{u})-C^{\#}(s,\vect{u})\right) + O(n^{-1/2}),$$  
where $C^{\#}(s,\vect{u})=sC(\vect{u})$, 
is  called sequential empirical process, see e.g. Chapter 2.12 in \cite{vandwell1996}. The remainder term $O(n^{-1/2})$, holds uniformly in $s, \vect{u}$ and over all copulas $C$. Estimating the marginal distribution functions in $\alpha_n^{\#}$ by their empirical counterparts we arrive at the so-called 
sequential empirical copula processes defined as
\begin{align*}
	\Cb_n^{\#}(s,\vect{u}) &= \frac{1}{\sqrt{n}}\sum_{i=1}^{\lfloor sn \rfloor} \left( \ind\{ \vect{X}_i \leq F_n^-(\vect{u})  \} - C(\vect{u}) \right)  =  \sqrt{n}\left( C_n^{\#}(s,\vect{u}) - \frac{\lfloor sn \rfloor}{n} C(\vect{u})\right) \\
	&= \sqrt{n}\left( C_n^{\#}(s,\vect{u}) - C^{\#}(s,\vect{u})\right) +O(n^{-1/2}) \\
	\Cb_n^+(s,\vect{u}) &=\frac{1}{\sqrt{n}}\sum_{i=1}^{\lfloor sn \rfloor} \left( \ind\{ \vect{X}_i \leq F_n^-(\vect{u})  \} - C_n(\vect{u}) \right) =  \sqrt{n}\left( C_n^{\#}(s,\vect{u}) - \frac{\lfloor sn \rfloor}{n} C_n^{\#}(1,\vect{u})\right) \\
	&= \sqrt{n}\left( C_n^{\#}(s,\vect{u}) - s C_n(\vect{u})\right) +O_\Prob(n^{-1/2})
\end{align*}
where $F_n^-(\vect{u})=(F_{n1}^-(u_1), \dots, F_{nd}^-(u_d))$ and $C_n^{\#}(s,\vect{u})  =n^{-1}\sum_{i=1}^{\lfloor sn \rfloor} \ind\{ \vect{X}_i \leq F_n^-(\vect{u})  \}$. For $s=1$ we retrieve the usual empirical copula process from Section \ref{sec:empcop}, i.e. $\Cb_n^{\#}(1,\cdot)=\Cb_n$. The process $\Cb_n^{\#}$ was introduced and investigated  in \cite{rueschendorf1976}. In order to state his general result we will make the following assumption on the weak convergence of the sequential empirical process $\alpha_n^{\#}$, which is similar to Condition \ref{cond:weak} in Section \ref{sec:empcop}.

\begin{cond}\label{cond:weakseq}
	For the strictly stationary sequence $(\vect{X}_i)_{i\in\Z}$  it holds that 
	$\alpha_n^{\#}=\sqrt{n}(G_n^{\#} - C^{\#})$ weakly converges in $\ell^\infty([0,1]^{d+1})$ to a tight, centered Gaussian field $\Bb_C^{\#}$ concentrated on $\Db_0^{\#}$, where
\begin{align*} 
	\Db_0^{\#}=\left\{ \alpha\in C([0,1]^{d+1})\, | \, \alpha(0,\cdot)=0 \text{ and } \alpha(s,\cdot) \in\Db_0 \text{ for all } s\in(0,1] \right\}.
\end{align*}
\end{cond}

As for the usual empirical process, the literature provides several sufficient conditions for the weak convergence postulated in Condition \ref{cond:weakseq}.
 For instance, in the case of i.i.d. random vectors we have weak convergence to the  Kiefer-M\"uller process, which has mean zero and covariance function
$$
	\Cov(\Bb_C^{\#}(s,\vect{u}), \Bb_C^{\#}(t,\vect{v})) = (s\wedge t) \left( C(\vect{u}\wedge \vect{v}) - C(\vect{u})C(\vect{v})\right),
$$
see e.g. Chapter 2.12 in \cite{vandwell1996}. For a similar result under strongly mixing conditions see e.g. \cite{philpinz1980} [which again incorporates the ARMA$(p,q)$-process discussed in Example \ref{ex:arma}.]

For technical reasons we consider the asymptotic behavior of $\Cb_n^+$ first. For $H^{\#}\in \ell^\infty([0,1]^{d+1})$ let $H(\vect{x})=H^{\#}(1,\vect{x})$ and consider the mapping 
\begin{align*} 
\Psi: \begin{cases}
	\Db_\Psi &\rightarrow \ \ \ell^\infty([0,1]^{d+1})\\
	H^{\#} &\mapsto \ \ \Psi(H^{\#})(s,\vect{x})=H^{\#}(s,H^-(\vect{x})) - s\, H( H^-(\vect{x}) )
\end{cases}
\end{align*}
where $\Db_\Psi$ consists of all functions $H^{\#}\in \ell^\infty([0,1]^{d+1})$ such that $H(\cdot)=H^{\#}(1,\cdot)\in\Db_\Phi$.
Recall from  Section \ref{sec:empcop} that $C_n(\vect{u})=G_n(G_n^-(\vect{u}))$. Similarly, it holds that $C_n^{\#}(s,\vect{u})=G_n^{\#}(s, G_n^-(\vect{u}))$ and thus 
\begin{align*}
\Cb_n^+(s,\vect{u} ) &= \sqrt{n}\left( G_n^{\#}(s, G_n^-(\vect{u})) - s\, G_n(G_n^-(\vect{u})) \right)  + O_\Prob(n^{-1/2}) \\
&= \sqrt{n} \left( \Psi(G_n^{\#}) - \Psi(C^{\#}) \right)(s,\vect{u})+ O_\Prob(n^{-1/2}),
\end{align*}
where the last equality follows from the fact that $\Psi(C^{\#})(s,\vect{u})=C^{\#}(s,\vect{u})-s\,C^{\#}(1,\vect{u})=0$. Similarly, define
\begin{align*}
\Gamma: \begin{cases}
	\Db_\Psi &\rightarrow \ \ \ell^\infty([0,1]^{d+1})\\
	H^{\#} &\mapsto \ \ \Gamma(H^{\#})(s,\vect{x})=H^{\#}(s,H^-(\vect{x})) 
\end{cases}
\end{align*}
and note that $\Cb_n^{\#} = \sqrt{n} \left( \Gamma(G_n^{\#}) - \Gamma(C^{\#}) \right) + O_\Prob(n^{-1/2})$. Mimicking the argumentation in Section \ref{sec:empcop}, in order to derive the asymptotics of $\Cb_n^+$ and $\Cb_n^{\#}$ it remains to show Hadamard-differentiability of $\Psi$ and $\Gamma$ at $C^{\#}$ [tangentially to suitable subspaces, if necessary]. This is done in the following Theorem. Note, that for the first part we do not need any assumptions on the smoothness of $C$, the result holds for every copula $C$.

\begin{theorem} \label{theo:Cseq_hadamard}
a) Let $C$ be an arbitrary copula.  Then the mapping $\Psi$ is Hadamard-differentiable at $C^{\#}\in\Db_\Psi$, $C^{\#}(s,\vect{x})=sC(s,\vect{x})$ with derivative
$$\Psi'_{C^{\#}}(\alpha^{\#})(s,\vect{x}) = \alpha^{\#}(s,\vect{x}) - s \alpha^{\#}(1,\vect{x}).$$
b) If moreover the copula satisfies condition \ref{cond:pd}, then the mapping $\Gamma$ is Hadamard-differentiable at $C^{\#}$ tangentially to $\Db_0^{\#}$. Its derivative at $C^{\#}$ in $\alpha^{\#}\in\Db_0^{\#}$ is given by
$$ \Gamma_{C^{\#}}'(\alpha^{\#})(s,\vect{x})=\alpha^{\#}(s,\vect{x}) - s \sum_{p=1}^d \partial_pC(\vect{x})\ \alpha^{\#}(1,\vect{x}^{(p)}). $$ 
\end{theorem}

\textbf{Proof. } 
We begin with the proof of \textit{a)}. Let $t_n\rightarrow0$ and $\alpha_n^{\#}\in \ell^\infty([0,1]^{d+1})$ with $\alpha_n^{\#}\rightarrow \alpha^{\#} \in \ell^\infty([0,1]^{d+1})$ such that $C^{\#}+t_n\alpha_n^{\#} \in \Db_\Psi$. Then
\begin{align*}
 &~ t_n^{-1} \{ \Psi(C^{\#}+t_n\alpha_n^{\#}) - \Psi(C^{\#}) \} (s,\vect{u}) \\
= &~ t_n^{-1} \{ (C^{\#}+t_n\alpha_n^{\#})(s,(C+t_n\alpha_n)^-(\vect{u}) )  - s(C+t_n\alpha_n)) ((C+t_n\alpha_n)^-(\vect{u})) \} \\
= &~ \alpha_n^{\#}(s, (C+t_n\alpha_n)^-(\vect{u}) ) - s\alpha_n((C+t_n\alpha_n)^-(\vect{u})).
\end{align*}
By the proof of Theorem \ref{theo:C_hadamard} we can conclude that $(C+t_n\alpha_n)^-$ uniformly converges to the identity on $[0,1]^d$, see equation \eqref{PsiU}, whose proof did not need Condition \ref{cond:weak}. Together with uniform convergence of $\alpha_n^{\#}$ and uniform continuity of $\alpha^{\#}$ we obtain assertion \textit{a)}. 

For the proof of \textit{b)} note that $\Gamma(\alpha^{\#})(s,\vect{x})=\Psi(\alpha^{\#})(s,\vect{x}) + s\, \Phi(\alpha)(\vect{x})$ with $\Psi$ and $\Phi$ from Theorem \ref{theo:Cseq_hadamard} \textit{a)} and Theorem \ref{theo:C_hadamard}, respectively. Some thoughts reveal that the Hadamard-differentiability of $\Psi$ and $\Phi$ transfers to $\Gamma$, such that $\Gamma'_{C^{\#}}(\alpha^{\#})(s,\vect{x})=\Psi'_{C^{\#}}(\alpha^{\#})(s,\vect{x}) + s\, \Phi'_{C}(\alpha)(\vect{x}).$
The details are omitted for the sake of brevity.  \qed

\medskip

The following main result of this section is a consequence of Theorem \ref{theo:Cseq_hadamard} and the functional delta method. Again note, that we do not need any smoothness assumptions on $C$ for the first part of the result.

\begin{corollary}\label{cor:empcopseq}
a) Suppose Condition \ref{cond:weakseq} holds. Then we have, for any copula $C$, that
\begin{align*}  
	\Cb_n^+(s,\vect{u}) =  \sqrt{n}\left( C_n^{\#}(s,\vect{u}) - s C_n(\vect{u})\right) +O_\Prob(n^{-1/2}) \weak \Gb_C^+    ~~~\text{ in }~~~ l^{\infty}([0,1]^{d+1}), 
\end{align*}
where the limiting Gaussian field may be expressed as
$\mathbb{G}_C^+(s,\vect{u})=\Bb_C^{\#}(s,\vect{u}) - s \Bb_C^{\#}(1,\vect{u}).$ 

b) If moreover the copula satisfies Condition \ref{cond:pd}, then we have
\begin{align*}  
	\Cb_n^{\#}(s,\vect{u}) =  \sqrt{n}\left( C_n^{\#}(s,\vect{u}) - s C(\vect{u})\right) +O_\Prob(n^{-1/2}) \weak \Gb_C^{\#}    ~~~\text{ in }~~~ l^{\infty}([0,1]^{d+1}), 
\end{align*}
where $\mathbb{G}_C^{\#}(s,\vect{u})=\Bb_C^{\#}(s,\vect{u}) - s \sum_{p=1}^d \partial_pC(\vect{u})\ \Bb_C^{\#}(1,\vect{u}^{(p)}).$ 
\end{corollary}

Corollary \ref{cor:empcopseq} has several applications for the investigation of multivariate  time series, when one is interested in tests for structural breaks in the dependence structure. For instance, \cite{wiedehvanvog2011} test for a constant Spearman's rho over time of some strictly stationary time series, while \cite{remillard2010} tests for a constant copula. The results of both papers are improved by Corollary \ref{cor:empcopseq} $a)$, since we do not make any smoothness assumption on $C$ at all.

\bigskip

{\bf Acknowledgements} 
The authors would like to thank \textit{Johan Segers} for fruitful discussions and for encouraging us to write this paper; and \textit{Martin Ruppert} for drawing our attention to the sequential empirical copula process $\Cb_n^+$ in more detail, especially without imposing any smoothness assumption at all.

This work has been supported   by the Collaborative
Research Center ``Statistical modeling of nonlinear dynamic processes'' (SFB 823, Teilprojekt A1, C1) of the German Research Foundation (DFG).

\bibliography{funcdelta}

\begin{thebibliography}{}

\bibitem[B{\"{u}}cher, 2011]{buecher2011}
B{\"{u}}cher, A. (2011).
\newblock {\em Statistical inference for copulas and extremes.}
\newblock PhD thesis, Ruhr-University Bochum, Germany.

\bibitem[B{\"{u}}cher and Dette, 2010]{buecdett2010}
B{\"{u}}cher, A. and Dette, H. (2010).
\newblock A note on bootstrap approximations for the empirical copula process.
\newblock {\em Statist. Probab. Lett.}, 80:1925--1932.

\bibitem[B\"{u}hlmann, 1993]{buehlmann1993}
B\"{u}hlmann, P.~L. (1993).
\newblock {\em The blockwise bootstrap in time series and empirical processes}.
\newblock ProQuest LLC, Ann Arbor, MI.
\newblock Thesis (Dr.Sc.Math)--Eidgenoessische Technische Hochschule Z\"{u}rich
  (Switzerland).

\bibitem[Dehling and Durieu, 2011]{dehlduri2011}
Dehling, H. and Durieu, O. (2011).
\newblock Empirical processes of multidimensional systems with multiple mixing
  properties.
\newblock {\em Stochastic Processes and their Applications}, 121:1076--1096.

\bibitem[Doukhan et~al., 2009]{douferlan2009}
Doukhan, P., Fermanian, J.-D., and Lang, G. (2009).
\newblock An empirical central limit theorem with applications to copulas under
  weak dependence.
\newblock {\em Stat. Inference Stoch. Process.}, 12(1):65--87.

\bibitem[Doukhan et~al., 1995]{doumasrio1995}
Doukhan, P., Massart, P., and Rio, E. (1995).
\newblock Invariance principles for absolutely regular empirical processes.
\newblock {\em Ann. Inst. H. Poincar\'e Probab. Statist.}, 31(2):393--427.

\bibitem[Durieu and Tusche, 2011]{duritusc2011}
Durieu, O. and Tusche, M. (2011).
\newblock An empirical process central limit theorem for multidimensional
  dependent data.
\newblock {\em arXiv:1110.0963}.

\bibitem[Fermanian et~al., 2004]{ferradweg2004}
Fermanian, J.-D., Radulovi{\'c}, D., and Wegkamp, M. (2004).
\newblock Weak convergence of empirical copula processes.
\newblock {\em Bernoulli}, 10(5):847--860.

\bibitem[Gaenssler and Stute, 1987]{gaenstut1987}
Gaenssler, P. and Stute, W. (1987).
\newblock {\em Seminar on empirical processes}, volume~9 of {\em DMV Seminar}.
\newblock Birkh\"auser Verlag, Basel.

\bibitem[Genest et~al., 2011]{gennesque2011}
Genest, C., Ne{\v{s}}lehov{\'a}, J., and Quessy, J.-F. (2011).
\newblock Tests of symmetry for bivariate copulas.
\newblock {\em The Annals of the Institute of Statistical Mathematics}, 64.
\newblock in press.

\bibitem[Kosorok, 2008]{kosorok2008}
Kosorok, M.~R. (2008).
\newblock {\em Introduction to Empirical Processes and Semiparametric
  Inference}.
\newblock Springer, New York.

\bibitem[Marinucci, 2005]{marinucci2005}
Marinucci, D. (2005).
\newblock The empirical process for bivariate sequences with long memory.
\newblock {\em Stat. Inference Stoch. Process.}, 8(2):205--223.

\bibitem[Mokkadem, 1988]{mokkadem1988}
Mokkadem, A. (1988).
\newblock Mixing properties of {ARMA} processes.
\newblock {\em Stochastic Process. Appl.}, 29(2):309--315.

\bibitem[Nelsen, 2006]{nelsen2006}
Nelsen, R.~B. (2006).
\newblock {\em An introduction to copulas}.
\newblock Springer Series in Statistics. Springer, New York, second edition.

\bibitem[Philipp and Pinzur, 1980]{philpinz1980}
Philipp, W. and Pinzur, L. (1980).
\newblock Almost sure approximation theorems for the multivariate empirical
  process.
\newblock {\em Z. Wahrsch. Verw. Gebiete}, 54(1):1--13.

\bibitem[Pollard, 1984]{pollard1984}
Pollard, D. (1984).
\newblock {\em Convergence of Stochastic Processes}.
\newblock Springer, New York.

\bibitem[R\'{e}millard, 2010]{remillard2010}
R\'{e}millard, B. (2010).
\newblock Goodness-of-fit tests for copulas of multivariate time series.
\newblock Technical report, HEC Montr\'{e}al.

\bibitem[R{\'e}millard and Scaillet, 2009]{remiscai2009}
R{\'e}millard, B. and Scaillet, O. (2009).
\newblock Testing for equality between two copulas.
\newblock {\em J. Multivariate Anal.}, 100(3):377--386.

\bibitem[Rio, 2000]{rio2000}
Rio, E. (2000).
\newblock {\em Th\'eorie asymptotique des processus al\'eatoires faiblement
  d\'ependants}, volume~31 of {\em Math\'ematiques \& Applications (Berlin)
  [Mathematics \& Applications]}.
\newblock Springer-Verlag, Berlin.

\bibitem[Ruppert, 2011]{ruppert2011}
Ruppert, M. (2011).
\newblock Consistent testing for a constant copula under strong mixing based on
  the tapered block multiplier technique.
\newblock Technical report.
\newblock unpublished.

\bibitem[R{\"{u}}schendorf, 1976]{rueschendorf1976}
R{\"{u}}schendorf, L. (1976).
\newblock Asymptotic distributions of multivariate rank order statistics.
\newblock {\em Annals of Statistics}, 4:912--923.

\bibitem[Schmid and Schmidt, 2007]{schmschm2007}
Schmid, F. and Schmidt, R. (2007).
\newblock Multivariate extensions of {S}pearman's rho and related statistics.
\newblock {\em Statist. Probab. Lett.}, 77(4):407--416.

\bibitem[Segers, 2011]{segers2011}
Segers, J. (2011).
\newblock Asymptotics of empirical copula processes under nonrestrictive
  smoothness assumptions.
\newblock {\em arXiv:1012.2133v}.

\bibitem[Sklar, 1959]{sklar1959}
Sklar, A. (1959).
\newblock Fonctions de r{\'{e}}partition {\`{a}} {$n$} dimensions et leurs
  marges.
\newblock {\em Publ.\ Inst.\ Statist.\ Univ.\ Paris}, 8:229--231.

\bibitem[Tsukahara, 2005]{tsukahara2005}
Tsukahara, H. (2005).
\newblock Semiparametric estimation in copula models.
\newblock {\em Canad. J. Statist.}, 33(3):357--375.

\bibitem[{V}an~der Vaart and Wellner, 1996]{vandwell1996}
{V}an~der Vaart, A.~W. and Wellner, J.~A. (1996).
\newblock {\em Weak Convergence and Empirical Processes - Springer Series in
  Statistics}.
\newblock Springer, New York.

\bibitem[van~der Vaart and Wellner, 2007]{vandwell2007}
van~der Vaart, A.~W. and Wellner, J.~A. (2007).
\newblock Empirical processes indexed by estimated functions.
\newblock In {\em Asymptotics: particles, processes and inverse problems},
  volume~55 of {\em IMS Lecture Notes Monogr. Ser.}, pages 234--252. Inst.
  Math. Statist., Beachwood, OH.

\bibitem[Wied et~al., 2011]{wiedehvanvog2011}
Wied, D., Dehling, H., van Kampen, M., and Vogel, D. (2011).
\newblock A fluctuation test for constant spearman’s rho.
\newblock Technical Report 16/2011, Technische Universit\"{a}t Dortmund, SFB
  823 Discussion Papers.

\end{thebibliography}

\end{document}